\documentclass[11pt]{amsart}
\usepackage{amssymb,latexsym}
\newtheorem{theorem}{Theorem}[section]
\newtheorem{prop}[theorem]{Proposition}
\newtheorem{lemma}[theorem]{Lemma}

\newtheorem{question}[theorem]{Question}
\newtheorem{definition}[theorem]{Definition}
\newtheorem{cor}[theorem]{Corollary}
\newtheorem{conjecture}[theorem]{Conjecture}
\newtheorem{speculation}[theorem]{Speculation}

\begin{document}

\title{The space of symplectic structures on closed 4-manifolds}

\author{Tian-Jun Li}
\address{School  of Mathematics\\  University of Minnesota\\ Minneapolis, MN 55455}
\email{tjli@math.umn.edu}

\thanks{The author is supported in part by NSF grant 0435099 and the McKnight fellowship.}
\maketitle

\section{Introduction}

Let $X$ be a $2n-$dimensional smooth manifold. A $2-$form $\omega$
on $X$ is said to be non-degenerate if, for each $q\in X$ and for
each nonzero vector $v$ in the tangent space $T_qX$, there is a
tangent vector $v\in T_qX$ such that $\omega(u,v)\ne 0$. A
symplectic structure on $X$ is a non-degenerate closed $2-$form. The
fundamental example of a symplectic structure is $\omega_0=\sum_i
dx_i\wedge dy_i$ on ${\bf R}^{2n}=\{(x_1,y_1, ...,x_n,y_n)\}$. In
fact, by the Darboux Theorem,
 every symplectic
structure  is locally like $({\bf R}^{2n},\omega_0)$.

Symplectic structures first appeared in Hamiltonian mechanics. A
K\"ahler form on a complex manifold is symplectic, thus we also find
a rich source of symplectic manifolds in algebraic geometry. Thirty
years ago people even wondered whether there are closed non-K\"ahler
symplectic manifolds. We have now gradually  realized that the world
of K\"ahler manifolds only occupies a tiny part of the symplectic
world (see \cite{Go1}, \cite{OT} and the references therein).

Two of the basic questions about symplectic structures are (see
\cite{McS2}): 1. Which smooth manifolds support symplectic
structures? 2. How many symplectic structures, up to appropriate
equivalence, are there on a given smooth manifold? In this survey we
focus on the second question for closed smooth $4-$manifolds (for
the first question in dimension $4$ see \cite{L4}).

In section 2 we review some fundamental  facts about symplectic
structures. In section 3 we survey what is known about the space of
symplectic structures in dimension 4. In the case of $b^+=1$ we have
a rather good understanding.  Especially, for a rational or ruled
manifold, there is the deep uniqueness result that a symplectic form
is determined by its cohomology class up to diffeomorphisms. We
further
 give a simple description of  the moduli spaces for such a manifold.
 We also point out various possible extensions to the case $b^+>1$.
 In section 4 we compare the space of symplectic
forms  and the  space of K\"ahler manifolds on a manifold admitting
a K\"ahler structure. As a by product we describe an example of
non-holomorphic Lefschetz fibration on a K\"ahler surface.

We would like to thank  P. Biran, J. Dorfmeister, R. Friedman and M.
Usher for useful discussions. We are very grateful to  the referee
for reading it carefully and making  many useful suggestions.

\section  {Symplectic structures}

In this section $X$ is a smooth, oriented, closed $2n-$dimensional
manifold.
\subsection
{Three spaces: $\Omega(X), {\mathcal M}_X, \bar{\Omega}(X)$}
  It is easy to see that  a $2-$form $\omega$ on $X$ is non-degenerate
  if and only if  $\omega^n$ is a volume form on $X$. Thus a non-degenerate form on
  $X$
  determines an orientation of $X$. If this orientation agrees with the given
one, we say that $\omega$ is orientation$-$compatible.

\begin{definition} Let ${\Omega}(X)$ denote the space of
orientation$-$compatible symplectic forms on $X$.
\end{definition}

$\Omega(X)$ is invariant under scaling by a positive scalar, i.e. if
$\omega \in {\Omega}(X)$, then ${\mathbb R}^+\omega\subset {\Omega}(X)$.
Moreover, ${\mathbb R}^-\omega \subset {\Omega}(X)$ if $n$ is even.

 Since the non-degeneracy condition is an open condition and $X$ is
compact, we see that, for any $\omega\in \Omega(X)$ and a
sufficiently small closed $2-$form $\eta$, the closed form
$\omega+\eta$ is still non-degenerate. Thus ${\Omega}(X)$ is an open
submanifold of the space of closed $2-$forms. In fact, we have

\begin{lemma} \label{convex} For every  $\omega\in \Omega(X)$, $\Omega(X)$ contains  a convex neighborhood
of $\omega$ in the space of closed $2-$forms.
\end{lemma}

\begin{proof}
Choose a Riemannian metric $g$ on $X$ and use it to identify $T^*X$
and $TX$. Now each closed $2-$form $\eta$ is viewed as an
endomorphism $L_{\eta}$ of  the inner product bundle $TX$ and $\eta$
is symplectic if and only if $L_{\eta}$ is an automorphism. Let
$\|L_{\eta}\|$ be the operator norm of $L_{\eta}$, which is defined
by
$$\|L_{\eta}\|=\hbox{min}_{|u|=1}|L_{\eta}(u)|.$$
Define also the radius of $L_{\eta}$ to be
$$r(L_{\eta})=\hbox{max}_{|u|=1}|L_{\eta}(u)|.$$
Since $X$ is compact both $\|L_{\eta}\|$ and $r(L_{\eta})$ are
finite. Moreover,
\begin{equation} \label{norm}
\begin{array}{ll}
\|L_{\eta_1+\eta_2}\|&\geq \|L_{\eta_1}\|-r(L_{\eta_2}),\cr
r(L_{t_1\eta_1+t_2\eta_2})&\leq
r(L_{t_1\eta_1})+r(L_{t_2\eta_2})=t_1r(L_{\eta_1})+t_2r(L_{\eta_2}),\cr
\end{array}
\end{equation}

 For each $\omega\in \Omega(X)$ let
  $U_{\omega}$  be the ball in the space of closed $2-$forms,
  $$U_{\omega}=\{\, \omega+\eta \,\, | \,\,  r(L_{\eta})< \frac {\|L_{\omega}\|}{2}\}.$$
   It follows from the second inequality of
(\ref{norm}) that for $\eta_1$ and $\eta_2$ in $U_{\omega}$, we have
$r(L_{t\eta_1+(1-t)\eta_2})< \frac {\|L_{\omega}\|}{2}$. Hence
$U_{\omega}$ is convex, and by the first inequality of (\ref{norm})
we have  $U_{\omega}\subset \Omega(X)$.
\end{proof}

Let  Diff$^+(X)$ be the group of orientation$-$preserving
diffeomorphisms of $X$. Then  $\Omega(X)$ is acted upon by
Diff$^+(X)$.

\begin{definition} \label{moduli}The quotient space ${\mathcal M}_X={\Omega}(X)/{\rm Diff}^+(X)$ is
called the moduli space of symplectic structures on $X$. Let $\bar
{\Omega}(X)$ be  the discrete set of  connected components of
${\mathcal M}_X$.
\end{definition}

As $\Omega(X)$ is an open subset of an infinite dimensional vector space,
a connected component is also path connected.
 Thus $\bar{\Omega}(X)$ is the same as the set of path connected components of
${\Omega}(X)$ up to the action of {\rm Diff}$^+(X)$. In some sense,
its cardinality is the number of symplectic structures on $X$.

Sometimes we will consider the reduced moduli space ${\mathcal
M}^r_X$, which is the quotient of the moduli space by the ${\mathbb
R}^+$ action.
\subsection {Cohomological invariants of symplectic structures}

Since $\omega$ is closed and $\omega^n$ is a volume form, $\omega$
represents a nonzero class in $H^{2}(X;{\mathbb R})$. In particular,
the second Betti number of $X$ is nonzero. Moreover,  in dimension
4, we must have $b^+\geq 1$ by the orientation-compatible
assumption.

Recall that an almost complex structure on $X$ is an endomorphism
$J$ on the tangent bundle $TX$ such that $J^2=-1$.  And an almost
complex structure $J$ is said to be compatible with $\omega$ if, for
any $q\in X$ and $u,v\ne 0\in T_qX$,
 $\omega(Ju, Jv)=\omega(u,v)$
and
 $\omega(u, Ju)>0$.
 The  space of $\omega-$compatible almost complex structures is non-empty and contractible.
Thus we can define the symplectic Chern classes
$c_i(X,\omega)=c_i(X,J)$, where $J$ is any $\omega-$compatible
almost complex structures.

\begin{definition}  $-c_1(X,\omega)\in H^2(X;{\mathbb Z})$ is called
the symplectic canonical class, and is denoted by $K_{\omega}$.
\end{definition}

Thus  there are two  basic cohomological invariants of $\omega$:
\[ \begin{array}{ll}
s:&{\Omega}(X)\longrightarrow H^2(X;{\mathbb R}),\quad
\omega\longrightarrow [\omega],\\
 c:&{\Omega}(X)\longrightarrow
H^2(X;{\mathbb Z}),\quad \omega\longrightarrow K_{\omega}.
\end{array} \]
Here $s$ stands for symplectic and $c$ stands for Chern.

\begin{definition} The image $s({\Omega}(X))\subset H^2(X;{\mathbb R})$ is called
the symplectic cone of $X$, and is denoted by ${\mathcal C}_X$. The
image $c({\Omega}(X))\subset H^2(X;{\mathbb Z})$ is denoted by
${\mathcal K}(X)$.
\end{definition}

Notice that ${\Omega}(X),{\mathcal M}_X, \bar\Omega(X), {\mathcal
C}_X, {\mathcal K}(X)$ are all differentiable invariants of $X$.
Notice also that, if $2n$ is divisible by $4$ and $K\in {\mathcal
K}(X)$, then $-K$ is also in ${\mathcal K}(X)$. This is because
$K_{-\omega}=-K_{\omega}$, and, in such a dimension, if $\omega$ is
orientation-compatible then $-\omega$ is also
orientation-compatible.

For each $\alpha\in {\mathcal C}_X$, let
$\Omega_{\alpha}(X)=s^{-1}(\alpha)$ be the space of symplectic
forms with cohomology class $\alpha$. For $\alpha\ne \beta$,
$\Omega_{\alpha}(X)$ is in general not homotopic to
$\Omega_{\beta}(X)$, so the map $s:\Omega(X)\rightarrow {\mathcal
C}_X$ is not a fibration. Nevertheless, non-trivial topology of
${\mathcal C}_X$ can sometimes be used to detect that of
$\Omega(X)$. A typical situation is the following: suppose
$\omega_u$ is a family of symplectic forms  parameterized by a
sphere and the associated family $[\omega_u]$  represents a nonzero
element of the homotopy groups of ${\mathcal C}_X$, then so does
$\omega_u$.

For each $K\in {\mathcal K}(X)$, let $\Omega(X,
K)=c^{-1}(K)=\{\omega\in \Omega(X)|K_{\omega}=K\}$ and ${\mathcal
C}_{(X, K)}=s(\Omega(X,K))$. Then $\Omega(X)$ is partitioned into the
disjoint union of $\Omega(X, K)$ over $K\in {\mathcal K}(X)$, and
likewise ${\mathcal C}_X$ is the union of ${\mathcal C}(X, K)$ over
$K\in {\mathcal K}(X)$ (not necessarily disjoint). ${\mathcal
C}_{(X,K)}$ is called the $K-$symplectic cone of $X$. It is an
analogue of the K\"ahler cones, but it is not clear that each of its
connected components is a convex cone (as the sum of two symplectic
forms may not be symplectic), though it is a union of rays.

Let D$(X)$ be the image of the group homomorphism
 $$\hbox{Diff}^+(X)\longrightarrow \hbox{Aut}(H^2(X;{\mathbb Z})),\quad \phi\longrightarrow \phi^*.$$
Since $K_{\phi^*\omega}=\phi^*K_{\omega}$ for any
$\phi\in$Diff$^+(X)$, D$(X)$ acts on ${\mathcal K}(X)$, and we have
a commutative diagram

\[\begin{array}{lllll}
 \hbox{Diff}^+(X)&\times& {\Omega}(X)&\longrightarrow &{\Omega}(X)\\
&\downarrow&&&\downarrow\\
  {\rm D}(X)&\times& {\mathcal
K}(X)&\longrightarrow &{\mathcal K}(X).  \end{array}
\]

\subsection{The equivalence relations}

There are several natural equivalence relations on $\Omega(X)$.

\begin{definition} Let $\omega_0$ and $\omega_1$ be in $\Omega(X)$. They are said to be
\begin{enumerate}
\item symplectomorphic if there exists a $\phi$ in {\rm Diff}$^+(X)$ such
that $\phi^*\omega_1=\omega_0$.

\item deformation equivalent if they can be joined by a smooth family
of symplectic forms $\omega_t$ on $X$.

\item isotopic, or strongly deformation equivalent,
 if they can be joined by a smooth family of
cohomologous symplectic forms $\omega_t$ on $X$.

\item strongly isotopic if there is a smooth isotopy $\psi_t$ of $X$ such
that $\psi_1^*\omega_1=\omega_0$. \footnote{Here we are following
\cite{McS2}. It also makes sense to call (3) weakly isotopic and (4)
isotopic as suggested by the referee.}

\end{enumerate}
\end{definition}

Clearly,  strongly isotopic forms are isotopic and symplectomorphic,
and isotopic forms are  cohomologous and  deformation equivalent. In
addition, deformation equivalent symplectic structures have the same
symplectic canonical class. Thus, from the commutative diagram in
the previous page, we have

\begin{lemma} \label{surjection} If  we let $\bar {\mathcal K}(X)$
denote the quotient of ${\mathcal K}(X)$ by {\rm D}$(X)$, then $\bar
{\Omega}(X)$ maps onto $\bar {\mathcal K}(X)$.
\end{lemma}

The following Moser stability is a fundamental result.

\begin{theorem} If $X$ is closed and $\omega_t, t\in [0,1]$, is a smooth family of cohomologous
symplectic forms on $X$, then there is a smooth isotopy $\psi_t$ of
$X$ such that $\psi_t^*\omega_t=\omega_0$.
\end{theorem}

It follows from the Moser stability  that isotopic forms are also
strongly isotopic when $X$ is closed (this is definitely not the
case when $X$ is not closed).
  If we let Diff$_0(X)$ denote the
group of diffeomorphisms isotopic to the identity, then the Moser
stability is equivalent to that Diff$_0(X)$ acts transitively on
each connected component (same as path connected here) of
$\Omega_{\alpha}(X)$. Thus, for each $\omega$, if
$\Omega_{[\omega]}(X, \omega)$ denote the connected component of
$\Omega_{[\omega]}$ containing $\omega$, then
$$\Omega_{[\omega]}(X,\omega)\cong \hbox{Diff}_0(X)/\hbox{Symp}(X,\omega)\cap \hbox{Diff}_0(X).$$
The associated homotopy long exact sequence to the corresponding bundle
\begin{equation}
\hbox{Symp}(X,\omega)\cap \hbox{Diff}_0(X)\longrightarrow
\hbox{Diff}_0(X)\longrightarrow \Omega_{[\omega]}(X,\omega)
\label{fibration}
\end{equation}
can then be used to analyze the topology of $\Omega_{[\omega]}(X,
\omega)$ in terms of that of Symp$(M,\omega)$, and vice versa.

\subsection{Manifold structure on $\mathcal M_X$}
Another  consequence of the Moser stability is the existence of a
manifold structure on ${\mathcal M}_X$. Let $\tilde {\mathcal
M}_X=\Omega(X)/$Diff$_0(X)$. Then $s$ descends to a map $\tilde
s:\tilde {\mathcal M}_X\to H^2(X;{\mathbb R})$.

 \begin{lemma} \label{local} $\tilde s$ is a local homeomorphism.
 \end{lemma}

 \begin{proof} For each $\omega\in \Omega(X)$ consider a convex open
 neighborhood $U_{\omega}$ constructed in Lemma \ref{convex}.
 Then $U'_{\omega}=$Diff$_0(X)\cdot U_{\omega}$ is an invariant open
 set of $\Omega(X)$
 for the action of Diff$_0(X)$. Let $\tilde U_{\omega}$ be the
 corresponding quotient, which is an open subset of $\tilde {\mathcal M}_X$.
 Since $U_{\omega}$ is convex all
 cohomologous symplectic forms in $U_{\omega}$ are on the same orbit
 of Diff$_0(X)$. Thus $\tilde s$ is a homeomorphism from
 $ \tilde U_{\omega}$ to the open set $s(U_{\omega})\subset H^2(X;{\mathbb
 R})$.

 \end{proof}

Let  $\Gamma(X)=$Diff$^+(X)/$Diff$_0(X)$ be the mapping class group
of $X$. Then $\mathcal M_X$ is  the discrete quotient $\tilde
{\mathcal M}_X/\Gamma(X)$.

\begin{cor} If $X$ is closed, then the moduli space ${\mathcal M}_X$ is
a (not necessarily Hausdorff or second countable) manifold of dimension $b_2(X)$.
\end{cor}

\subsection{Dimensions other than $4$}
For a closed oriented $2-$manifold $\Sigma$, $\bar {\Omega}(\Sigma)$
has only one point. Moreover, two cohomologous symplectic forms
which are deformation equivalent are isotopic, and hence strongly
isotopic. Therefore ${\mathcal M}_{\Sigma}$ is homeomorphic  to
${\mathbb R}^+$.

In dimensions  higher than $4$, pseudo-holomorphic curves techniques
introduced in \cite{Gr} have been used to distinguish symplectic
structures.

  Ruan
(\cite{R1}) constructed examples of manifolds in dimension 6 and
above with $\bar {\Omega}$ infinite by demonstrating that $|\bar
{\mathcal K}|$ being infinite. There are also cohomologous forms
with distinct symplectic canonical classes in dimension 6 and above
(\cite{R3}).

McDuff \cite{Mc4}  constructed examples of cohomologous symplectic
forms which are deformation equivalent but not isotopic in dimension
6 and above. There are even such forms which are not
symplectomorphic in dimension 8 and above.

 Suppose that $\omega_t$ and $\omega_t'$ are two families of
 symplectic forms on $X$ such that $\omega_0=\omega_0'$ and $\omega_t$ is cohomologous to $\omega_t'$
 for all $t\in [0,1]$. Then by the existence of a convex neighborhood
 of $\omega_0=\omega_0'$ there is some $\epsilon>0$ and an isotopy $\psi_t$ such that $\psi_t^*\omega_t'=\omega_t$
 for $0\leq t\leq \epsilon$. However, for the examples of $X$ in
 \cite{Mc4} there are $\omega_t$ and $\omega_t'$ where one can
 take $\epsilon$ to be any number less than $1$ but not equal to $1$.
 In fact, $\omega_1$ and $\omega_1'$ are not even symplectomorphic.
   In particular, for such an $X$, the moduli space ${\mathcal M}_X$ is not
Hausdorff as pointed out by the referee.

There are no such examples known in dimension 4.
\bigskip

\section{Space of symplectic structures on 4-manifolds}

In this section $X$ is a smooth, closed, oriented $4-$manifold with
$\Omega(X)$ non-empty. We would like to emphasize that we do not fix
a symplectic structure on $X$. Instead we are interested in the
space $\Omega(X)$. From this point of view,
 we need to introduce the notions of minimality and  the Kodaira dimension
for the manifold $X$ itself (rather for a pair $(X,\omega)$).

\subsection{Minimality and the Kodaira dimension $\kappa$}

 \begin{definition} Let ${\mathcal E}_X$ be the set of cohomology classes whose
Poincar\'e dual are represented by smoothly embedded spheres of
self-intersection $-1$. $X$ is said to be (smoothly) minimal if
${\mathcal E}_X$ is the empty set.
 \end{definition}

Equivalently, $X$ is minimal if it is not the connected sum of
another manifold  with $\overline{\mathbb {CP}}^2$. We say that $Y$
is  minimal model of $X$ if $Y$ is a minimal and $X$ is the
connected sum of $Y$ and a number of $\overline{\mathbb {CP}}^2$.

We  also recall the notion of minimality for a pair $(X,\omega)$
with $\omega\in \Omega(X)$. $(X,\omega)$  is said to be
(symplectically) minimal if ${\mathcal E}_{\omega}$ is the empty
set, where
$${\mathcal E}_{\omega}=\{E\in {\mathcal
E}_X|\hbox{ $E$ is represented by an embedded $\omega-$symplectic
sphere}\}.$$ Here is a relevant  and important fact:

\begin{theorem} \label{empty} If $Y=Z\# \overline{\mathbb {CP}}^2$, then $\Omega(Y)$ is
non-empty if and only if $\Omega(Z)$ is non-empty.
\end{theorem}
The if part is the easier part: it follows from the symplectic
blow-up construction. The only if part relies on the symplectic
blow-down construction and the following result (see \cite{T1},
\cite{LL1}, \cite{L3}): For any $\omega\in \Omega(X)$, ${\mathcal
E}_{\omega}$ is empty if and only if ${\mathcal E}_X$ is empty.

A rational $4-$manifold is $S^2\times S^2$ or ${\mathbb {CP}}^2\#
k{\overline{\mathbb {CP}}^2}$ for some non-negative integer $k$. A
ruled $4-$manifold is the connected sum of a number of (possibly
zero) ${\overline {\mathbb {CP}}^2}$ with an $S^2-$bundle over a
Riemann surface. A rational or a ruled manifold admits K\"ahler
structure and hence symplectic structures.

\begin{definition} \label{kodairadimension} $X$ is said to have Kodaira dimension
$\kappa=-\infty$ if $X$ is rational or ruled. \\
Otherwise, first suppose  $X$ is minimal. Then the Kodaira dimension
$\kappa$ of $X$ is defined in terms of ${\mathcal K}(X)$ as follows:
\begin{itemize}
\item
$\kappa(X)=0$ if there is a torsion $K\in {\mathcal K}(X)$;

\item  $\kappa(X)=1$ if there is a non-torsion $K\in {\mathcal K}(X)$
with $K^2=0$;

\item $\kappa(X)=2$ if there is a $K\in {\mathcal K}(X)$ with $K^2>0$.

\end{itemize}
For a general $X$, $\kappa(X)$ is defined to be $\kappa(Y)$, where
$Y$ is a minimal model of $X$ (this makes sense by Theorem
\ref{empty}).
\end{definition}

For the notion of the Kodaira dimension $\kappa(X,\omega)$ for a
pair $(X,\omega)$, see \cite{Le2}, \cite{McS}, \cite{L1}. It is
shown in \cite{L1} that $\kappa(X,\omega)$ only depends on the
oriented diffeomorphism type of $X$. This fact implies that
$\kappa(X)$ in Definition \ref{kodairadimension} is independent of
the choice of $K\in {\mathcal K}(X)$, and hence well-defined.

 Finally we mention that it is also useful to
consider, for $K\in {\mathcal K}(X)$, the following subset of
${\mathcal E}_X$,
$${\mathcal E}_{(X,K)}=\{E\in {\mathcal E}_X|E\cdot K=-1\}.$$

\subsection{Finiteness of $\bar{\Omega}(X)$ and $\bar{\mathcal
K}(X)$}

\begin{question} \label{finiteness} Is  $\bar {\Omega}(X)$ finite for every $X$?
\end{question}

This question asks whether there are only finite number of
symplectic structures up to deformations and diffeomorphisms. For
complex structures, the corresponding finiteness is true (see
\cite{FM2}). We can answer this question affirmatively only for
manifolds with $\kappa=-\infty$.

\begin{theorem} \label{uniqueness} If $\kappa(X)=-\infty$, then
$\bar{\Omega}(X)$ has only one element, i.e. {\rm Diff}$^+(X)$ acts
transitively on the connected components of $\Omega(X)$.
\end{theorem}

This  is a consequence of several deep results.  One key point why
such a classification is possible is the existence of embedded
spheres with non-negative self-intersections when
$\kappa(X)=-\infty$. See \cite{LMc1} for an excellent account in the
case of ${\mathbb CP}^2$ and $S^2-$bundles. For the non-minimal
ones, the additional inputs are Theorems \ref{canonicalclass} and
\ref{empty}.

In view of Lemma \ref{surjection}, the weaker question whether $\bar
{\mathcal K}(X)$ is finite is important regarding Question
\ref{finiteness}.
 This question was first answered when $b^+>1$. In this case, it is shown in
 \cite{T1} that every $K\in{\mathcal K}(X)$ is a Seiberg-Witten
 basic class. Together with the finiteness of  Seiberg-Witten
 basic classes in \cite{Wit}, we have

\begin{theorem}  If $b^+(X)>1$, then
 ${\mathcal K}(X)$ is finite.
 \end{theorem}

The case of $b^+=1$ was later handled in \cite{LL3}.
\begin{theorem}  \label{canonicalclass}  Let $X$ be a 4-manifold with $b^+=1$. Then $\bar {\mathcal K}(X)=
{\mathcal K}(X)/{\rm D}(X)$ has at most two elements. And if $X$ has
$\kappa=-\infty$, then $\bar {\mathcal K}(X)$ is a one point set.
 Moreover, if $X$ is minimal, then ${\mathcal K}(X)$ has at most two elements.
\end{theorem}

As an immediate consequence of the two theorems above, we have,

\begin{cor} $\bar {\mathcal K}(X)={\mathcal K}(X)/{\rm D}(X)$ is finite every any $X$.
\end{cor}

Some remarks are in order now.  The finiteness actually does not
hold for ${\mathcal K}(X)$: If $\kappa(X)=-\infty$ and $b^-(X)\geq
2$, then ${\mathcal K}(X)$ is an infinite set. If $\kappa(X)=0$,
then $\bar {\mathcal K}(X)$ is a one point set, and if $X$ is also
minimal, i.e. there is a torsion element $K$ in ${\mathcal K}(X)$,
then ${\mathcal K}(X)$ itself is a one point set
 (see \cite{McS},
\cite{L1}). In the case $b^+=1$, there are examples of minimal
manifolds with $\bar {\mathcal K}(X)={\mathcal K}(X)$ having two
elements (see \cite{KK2} for such an example due to Mostow
rigidity). When $b^+>1$, McMullen and Taubes (\cite{MT}) constructed
$X$ with  $\bar {\mathcal K}(X)/\pm
>1$, i.e. more than one symplectic canonical classes up to sign and
diffeomorphisms (see \cite{Sm1}, \cite{Le}, \cite{V1} for more such
examples).

Finally  we state a speculation,

\begin{speculation} $\bar {\mathcal K}(X)$ has
at most $2^{b^+}$ elements.
\end{speculation}

 For a $4-$manifold $X$ with no
embedded tori of self-intersection zero, this was raised in
\cite{FS}. For the general case see the discussion in \cite{L4}.
\subsection{The symplectic cone  ${\mathcal C}_X$ and the $K-$symplectic cone ${\mathcal C}_{(X,K)}$}

\begin{definition}The positive cone of $X$ is
$${\mathcal P}_{X}=\{e\in H^2(X;{\mathbb R})|e^2>0\}.$$
\end{definition}
By the orientation-compatible condition, the symplectic cone
${\mathcal C}_X$ is contained in the positive cone ${\mathcal P}_X$.
It is easy to see that the positive cone is  $(b^+-2)-$connected.
And $\pi_{b^+-1}({\mathcal P}_{X})={\mathbb Z}$ if $b^+>1$,
$\pi_{b^+-1}({\mathcal P}_{X})={\mathbb Z}_2$ if $b^+=1$ (see e.g.
\cite{LL2}).  A winding family of symplectic forms is an $S^{b^+-1}$
family of symplectic forms which represents a generator of
$\pi_{b^+-1}({\mathcal P}_X)$ (\cite{LL2}). As previously remarked,
such a family would also represent a nonzero element in
$\pi_{b^+-1}(\Omega(X))$.  Every known manifold of $\kappa=0$
carries such a family (\cite{L1}).

 For manifolds with $b^+=1$, the symplectic cone has been
completely determined in \cite{LL3}. Notice that in this case
${\mathcal P}_X$ has two path connected components. We start with
the $K-$symplectic cones.

\begin{theorem} \label{positivecone} Let $X$ be a 4-manifold with $b^+=1$ and $K\in {\mathcal K}(X)$.
\begin{enumerate}
\item  If $K$ is a torsion class, then ${\mathcal C}_{(X, K)}={\mathcal P}_X$.

\item  If $K$ is not a torsion class, then ${\mathcal C}_{(X,K)}$ is contained in
one of the components of ${\mathcal P}_X$, denoted by ${\mathcal
P}_{(X,K)}$. Moreover,
$${\mathcal C}_{(X,K)}=\{e\in {\mathcal P}_{(X,K)}|0<e\cdot E\hbox { for all
$E\in {\mathcal E}_{(X,K)}$}\}.$$

\item  In particular, ${\mathcal C}_{(X, K)}$ is one connected component
of the positive cone if $X$ is minimal and $K$ is not a torsion
class.

\item  Moreover, if $\Omega_{\omega}(X,\omega)$ is the connected component of
$\Omega(X)$ containing $\omega$, then
$s(\Omega_{\omega}(X,\omega))={\mathcal C}_{(X,K_{\omega})}$.

\item  If $K\ne K'$, then ${\mathcal C}_{(X, K)}$ and ${\mathcal
C}_{(X,K')}$ are disjoint.
\end{enumerate}
\end{theorem}

For the so called class ${\mathcal C}$ manifolds including
$4-$manifolds with $b^+=1$, $b_1=0$ and the ruled manifolds,  a
slightly weaker version of Theorem \ref{positivecone} appeared in
\cite{Bi3} and \cite{Mc3} (see also \cite{KK} for some partial
results on $S^1-$equivariant symplectic cones on rational
manifolds). One main technique in proving Theorem \ref{positivecone}
is the following inflation process first appeared in \cite{La} (see
also \cite{LMc2}): if $\omega\in \Omega(X)$ and $C\subset(X,\omega)$
is an embedded symplectic surface with $C\cdot C\geq 0$, then there
exists a closed $2-$form $\rho$ such that for every $t>0$, the form
$\omega_t=\omega+t\rho$ is in $\Omega(X)$. In the case $b^+=1$, the
inflation is very effective as there are always infinitely many
embedded symplectic surfaces from Taubes' SW$\Rightarrow$Gr,
together with the SW wall crossing formula in \cite{KM}, \cite{LL4},
\cite{OO}, \cite{OT2}. We remark that the idea to use the inflation
procedure in order to construct `new' symplectic forms on
4-manifolds by using symplectic curves  first appeared in \cite{Bi3}
and \cite{Bi4}.

For the symplectic cone itself, we have

\begin{theorem}  \label{symplecticcone} If $X$ has $b^+=1$,
then
$${\mathcal C}_X=\{e\in {\mathcal P_X}|0<|e\cdot E|\hbox { for all $E\in {\mathcal E}_X$}\}.$$
In particular, ${\mathcal C}_X$ is the positive cone if $X$ is
minimal.
\end{theorem}

Thus the  explicit description of   ${\mathcal C}_X$ is reduced to
the knowledge of ${\mathcal E}_X$.  In general  ${\mathcal E}_X$ is
easy to describe, and the difficult case of a rational manifold
 has been solved  in
\cite{LiL1}.

For manifolds with $b^+>1$, we have the following constraint (\cite{T1}).

\begin{theorem} \label{constraint} Let $X$ be a $4-$manifold with $b^+>1$ and $\omega\in \Omega(X)$.
Then  $\pm K_{\omega}$ is a SW basic class, and for any SW basic
class $B\ne -K_{\omega}\in H^2(X;{\mathbb Z})$, $[\omega]\cdot
(K_{\omega}+B)>0$. Therefore, for each $K\in {\mathcal K}(X)$,
\begin{equation} {\mathcal C}_{(X, K)}\subset \{e\in {\mathcal P}_X|e\cdot
(K+B)>0 \hbox{ for any SW basic class $B\ne -K$}\} \label{SWbasic}.
\end{equation}
\end{theorem}

It follows from (\ref{SWbasic}) and  part 5 of Theorem
\ref{positivecone} that we have
\begin{cor} For any $X$ and  $K\ne
K'\in {\mathcal K}(X)$, ${\mathcal C}_{(X, K)}$ and ${\mathcal
C}_{(X,K')}$ are disjoint.
\end{cor}

In \cite{LL3} we speculated that (\ref{SWbasic}) is an equality. A
 counterexample was constructed in \cite{V2}. However it is still possible that
(\ref{SWbasic}) is an equality with SW basic classes replaced by SW
monopole classes, as Theorem \ref{constraint} is still valid for  SW
monopole classes and the new obstruction in \cite{V2} may actually
come from a SW monopole class.
 Here a class $B\in H^2(X;{\mathbb Z})$ is said to be a SW monopole class if the SW moduli space
for a corresponding Spin$^c$ structure is non-empty for any pair of
a metric and a $2-$form on $X$.

 There are very few  cases where the
symplectic cone has been completely determined when $b^+>1$.
 For
orientable torus bundles over torus, it is shown in \cite{Ge} that
the symplectic cone is the positive cone. For K\"ahler surfaces see
the discussion in section 4.

If we can answer affirmatively the following question, then we would
be able to determine the symplectic cone of a manifold in terms of
its minimal model.

\begin{question} Suppose $X$ has $b^+>1$, $\beta\in{\mathcal C}_X$
and $e$ is a generator of $H^2(\mathbb{CP}^2;\mathbb{Z})$.
\begin{enumerate}
\item  Is $\beta-\lambda e$ in ${\mathcal C}_{X\#\overline{\mathbb{CP}}^2}$ for every
nonzero scalar $\lambda$ with $\lambda^2<\beta^2$?

\item  More generally, if $e_1, ...,e_k$ are generators of $H^2(\mathbb{CP}^2;\mathbb{Z})$
for $k$ copies of $\mathbb{CP}^2$, is
$\beta-\sum_{i=1}^k\lambda_ie_i$ in ${\mathcal
C}_{X\#k\overline{\mathbb{CP}}^2}$ for every $\{\lambda_i\}$ with
$\lambda_i\ne 0$ and $\sum_i\lambda_i^2<\beta^2$?
\end{enumerate}
\end{question}

Notice that, in the case $b^+=1$ and  $\kappa(X)\ne -\infty$, an
affirmative answer is provided by Theorem \ref{symplecticcone} (see
also \cite{Bi4}).

This question is closely related to the symplectic packing problem
(see \cite {McP}). In fact, if there is an $\omega\in \Omega(X)$
such that $\beta=[\omega]$ and $(X,\omega)$ has a full packing by
one standard symplectic ball, then part 1 is true. We can also
answer part 1 in the following situation:
 $\beta=[\omega]$ integral and $\beta^2=1$, and there is an embedded surface $C\subset X$,
 symplectic respect to $\omega$ and   representing the Poincar\'e dual of
$[\omega]$. To see this, consider a small blow up around a point in
$C$ to get a symplectic form $\tilde \omega$ on
$X\#\overline{\mathbb{CP}}^2$. The point is that we can construct
the proper transform of the $\omega-$symplectic surface  $C$, which
is an $\tilde \omega-$symplectic surface $\tilde C$ Poincar\'e dual
to $\beta-e$ (see \cite{MW2}). Since $\tilde C$ has
self-intersection zero, we can apply the inflation to get symplectic
forms whose classes are arbitrarily close to PD$(\tilde C)= \tilde
\omega-e$. Since $[\tilde \omega]$ can be chosen to be arbitrarily
close to $\beta=[\omega]$, we are done. This is just a special case
of the
 the following more general result, which is a consequence of the main result in \cite{Bi2} (see also
 \cite{Bi1})
 proved
using the symplectic packing of symplectic ruled surfaces.

\begin{theorem} \label{packing} Suppose $\beta$  is an integral  class in ${\mathcal C}_X$ and is represented by
an embedded $\omega-$symplectic surface $C$ of genus at least $1$
for some $\omega\in \Omega(X)$ with  $[\omega]=\beta$. Then, for any
integer $N\geq \beta^2$ and any nonzero scalar $\lambda$ with
$\lambda^2<\beta^2/k$, the class
$$\beta_{k,\lambda}=\beta-\lambda e_1-...-\lambda e_k$$ is in ${\mathcal
C}_{X\# k\overline{\mathbb {CP}}^2}$.
\end{theorem}

This theorem can also be seen as an application of the normal
connected sum construction in \cite{Go1}, \cite{MW}: summing $X$
with a non-minimal ruled manifold. In fact the same is true for the
inflation process. We will encounter another such application to
constructing  symplectic forms in the next section.

\subsection{Cohomologous symplectic forms}

Having determined the image of the map $s:\Omega_X\longrightarrow
H^2(X;{\mathbb R})$ in the case $b^+=1$, we are  also able to say
something about its inverse image (\cite{Mc3}, \cite{LL3}) by
another application of the inflation process.

\begin{theorem} \label{defiso}   Let $X$ be a $4-$manifold with
$b^+=1$, and $\omega_1$, $\omega_2$ be two cohomologous symplectic
forms on $X$. If $\omega_1$ and $\omega_2$ can be joined by a path
of symplectic forms, then they can be joined by a path of
cohomologous symplectic forms, i.e. if they are deformation
equivalent, then they are isotopic.
\end{theorem}

The case of $X=\mathbb CP^2\#l\overline{\mathbb CP}^2, l\leq 6$ was
essentially proved in \cite{Bi5}.  In terms of the map $s$, the
result above can be interpreted as saying that, when restricted to a
path connected component of $\Omega_X$, the inverse image of a point
of $s$ is path connected.

Recall that,  if $\Omega_{[\omega]}(X, \omega)$ denotes the
connected component of $\Omega_{[\omega]}$ containing $\omega$, then
$$\Omega_{[\omega]}(X,\omega)\cong \hbox{Diff}_0(X)/\hbox{Symp}(X,\omega)\cap \hbox{Diff}_0(X).$$
Seidel (\cite{Se}) has shown that there are many symplectic
manifolds $(X,\omega)$ such that Symp$(X,\omega)\cap$Diff$_0(X)$ is
not connected. It then follows from the  associated homotopy exact
sequence of the corresponding fibration that $\pi_1(\Omega_{[\omega]}(X,\omega))$ is non-trivial.
 Kronheimer (\cite{K}) further constructed directly higher dimensional
non-trivial families of cohomologous symplectic forms arising form
holomorphic rational curves with negative self-intersections. More
precisely,  for each positive integer $n\geq 2$, a holomorphic
rational curve with self-intersection $-n$ in an algebraic surface
sometimes gives rise to an $S^{2n-3}-$family of cohomologous
symplectic forms representing a non-trivial homology class in degree
$2n-3$. For $n=1$, such an example arises from a $2-$disc family of
quartics $X_t$ in ${\mathbb {CP}}^3$ with $X_t$ smooth when $t\ne 0$
and $X_0$ singular with a single ordinary double-point.   The
boundary $S^1-$family of smooth quartics is smoothly trivial as the
minimal resolution of the nodal quartic is again the K3 surface
(with an exceptional holomorphic rational $-2$ curve). The
non-trivial $S^1-$family of symplectic forms on the  K3 surface then
is simply the restriction of the Fubini-Study form on
$\mathbb{CP}^3$.
  We remark that  all the families of symplectic forms in \cite{K}
  are null-homologous in $\Omega(X)$.

 By Theorem \ref{uniqueness} and
Theorem \ref{defiso}, we have

\begin{theorem} \label{cohdiff}For $X$ with $\kappa=-\infty$,
cohomologous symplectic forms are symplectomorphic.
\end{theorem}

 Thus, for such an $X$, if Diff$(X,[\omega])$ is the subgroup of Diff$^+(X)$ preserving the
class $[\omega]$, then
$$\Omega_{[\omega]}(X)\cong \hbox{Diff}(X,[\omega])/\hbox{Symp}(X,\omega).$$
 The group Symp$(X,\omega)$ has been extensively
studied in \cite{Ab}, \cite{AM}. In particular, consider the case
$X=S^2\times S^2$, and $\omega_{\lambda}$ a product form
$(1+\lambda)\sigma_0\times \sigma_0$ with $\lambda>0$ \footnote{In
this case Symp$(X,\omega_{\lambda})$ is shown to be connected in
\cite{AM}, thus agreeing with
Symp$(X,\omega_{\lambda})\cap$Diff$^0(X)$. Notice also that in this
case Diff$(X,[\omega])$ agrees with the subgroup of Diff$^+(X)$
acting trivially on the homology of $X$.}. For an odd integer $l\geq
1$, the associated homotopy long exact sequence to the corresponding
fibre bundle is used in  \cite{AM} to identify a generator of the
$(2l-2)-$th homotopy group of Symp$(X,\omega_l)$ with Kronheimer's
non-trivial $(2l-1)-$dimensional family in $\Omega_{[\omega_l]}(X)$
associated to the unique $-(l+1)$ rational curve in the
Hirzebruch\footnote{$F_{l+1}$ is diffeomorphic to $S^2\times S^2$,
and $\omega_{l}$ is a K\"ahler form on $F_{l+1}$.} surface
$F_{l+1}$. There is a similar picture for the non-trivial
$S^2-$bundle $\mathbb {CP}^2\#\overline{\mathbb{CP}}^2$.

We would like to  point out that Theorem \ref{cohdiff} may play an
important role understanding symplectic 6--manifolds (\cite{LR}). We
give two examples here.  The first is a simple characterization of a
symplectic blowup at a point: A symplectic $6-$manifold $(N,\omega)$
is a symplectic blowup of another symplectic 6--manifold at a point
if and only if it contains an embedded symplectic 4--manifold
$(X,\eta)$ such that $X$ is diffeomorphic to $\mathbb{CP}^2$ and has
normal bundle ${\mathcal O}(-1)$. The second  provides a way to
recognize a uniruled 6--manifold: If a symplectic 6--manifold
$(N,\omega)$ contains an embedded symplectic 4--manifold $(X,\eta)$
such that $\kappa(X)=-\infty$ and $X$ has trivial normal bundle,
then $(N,\omega)$ is uniruled. Here $(N,\omega)$ is said to be
uniruled if there is a non-trivial one point genus zero GW invariant
\footnote{Notice that, e.g. by \cite{L2}, it is not meaningful to
define this notion by requiring that there is a symplectic sphere in
a fixed class through every point, otherwise every simply connected
manifold would be uniruled.}.

\subsection{$\Omega(X)$ and the  moduli space ${\mathcal M}_X$}

We have the following simple description of ${\mathcal M}_X$ in the
case of $b^+=1$.

\begin{theorem} \label{modulispace} Let $X$ be a 4-manifold with $b^+=1$. Then
${\mathcal M}_X$ consists of  ${|\bar {\Omega}(X)| \over {|\overline
{\mathcal K}(X)|}}$ number of connected components, each
homeomorphic to ${\mathcal C}_X/{\rm D}(X)$. In particular, if $X$
is minimal, then ${\mathcal M}_X$ consists ${ |\bar
{\Omega}(X)|\over {|\overline{\mathcal K}(X)|}}$ number of
${\mathcal P}_X/{\rm D}(X)$.
\end{theorem}

\begin{proof}  Pick a connected component $\Omega^0(X)$ of $\Omega(X)$ with symplectic canonical class $K$. We first consider the action
of Diff$_0(X)$ on $\Omega^0(X)$. By Theorem \ref{defiso}, each
inverse image of $s:\Omega^0(X)\rightarrow s(\Omega^0(X))$ is an
orbit of Diff$_0(X)$.
 Together with part 4 of Theorem \ref{positivecone}, we have
$$\Omega^0(X)/\hbox{Diff}_0(X)\cong s(\Omega^0(X))={\mathcal C}_{(X,K)}.$$

Now let $E\Omega(X)$ be the orbit of $\Omega^0(X)$ under the action of Diff$^+(X)$.
There are
$|\bar {\Omega}(X)|$ number of such orbits by the definition of $\bar \Omega (X)$.

Let us first assume that $\bar {\mathcal K}(X)$ is a one point set.
It suffices to show that  $E\Omega(X)$ gives rise to a copy of
${\mathcal C}_X/{\rm D}(X)$. Since Diff$^+(X)$ acts transitively on
${\mathcal K}(X)$, we have, by the preceding discussion,
$$E\Omega(X)/\hbox{Diff}_0(X)=s(E\Omega(X))=\cup_{K\in {\mathcal K}(X)}{\mathcal C}_{(X,K)}={\mathcal C}_X.$$
Now we are done with this case as the quotient group
Diff$^+(X)/\hbox{Diff}_0(X)$ acts on ${\mathcal C}_X \subset
H^2(X;{\mathbb R})$ via D$(X)$.

By Theorem \ref{canonicalclass}, the remaining case is when there are two elements
in ${\mathcal K}(X)$ represented by $K$ and $-K$. If $\omega$ is in an orbit
 $E\Omega(X)$, then there is a distinct orbit $E'\Omega(X)$ containing $-\omega$.
The argument above then shows that,
$$(E\Omega(X)\cup E'\Omega(X))/\hbox{Diff}_0(X)={\mathcal C}_X.$$

The last statement now follows from part 3 of Theorem
\ref{positivecone}.
\end{proof}

In particular, by Theorem \ref{uniqueness}, we have

\begin{theorem}
If  $\kappa(X)=-\infty$, then ${\mathcal M}_X$ is homeomorphic to
${\mathcal C}_X/{\rm D}(X)$.
\end{theorem}

In this case D$(X)$ has been explicitly determined in \cite{FM1},
\cite{LiL1} and \cite{LiL2} \footnote{It is also worth to mention
that D$(X)$ is finite when $X$ is of $b^+=1$ and Kodaira dimension 2
(see \cite{LL3}).}.

For $\mathbb {CP}^2$, the reduced moduli space is just a point.
For a
minimal ruled surface other than $S^2\times S^2$, the reduced moduli space is an open interval, while for $S^2\times S^2$ it is a
half-open interval because of the diffeomorphism interchanging the factors.
A natural question is whether there is a partial geometric
compactification of ${\mathcal M}^r_X$. We expect that  it is given
by forms degenerate along spheres of self-intersections $0$ or $1$
with canonical local models.

Finally we would like to mention two instances where certain
`geometric' symplectic forms can be determined up to isotopy.

Firstly, let $\pi:X\longrightarrow B$ be a fibration with $B$ an
oriented surface. $\omega\in \Omega(X)$ is said to be
$\pi-$compatible if all the fibers are $\omega-$symplectic surfaces.
It is shown in \cite{Mc3} all such forms are deformation equivalent.
Furthermore, when $X=F\times B$ with $g(F)\geq 1$, if $\omega$ is
$\pi-$compatible and cohomologous to a product form with the ratio
between the $\omega-$area of $B$ and the $\omega-$area of $F$ being
at most ${g(B)-1 \over {g(F)-1}}$, then $\omega$ is isotopic to a
split form.

Secondly, if $X$ has a free circle action, then it is shown in
\cite{HW} that the space of invariant symplectic forms is homotopic
equivalent to a subspace of non-degenerate $1-$forms on the quotient
$3-$manifold \footnote{By \cite{FGM}, $X/S^1$ is a fibred
$3-$manifold.}. A particular interesting case is, for the linear
action on the $4-$torus, this space is connected and simply
connected, and every invariant form is isotopic to a constant
coefficient form. Moreover, for certain cohomology classes $\alpha$,
$\Omega_{\alpha}(X)$ is shown to be homotopic to a circle.


\section{Symplectic forms versus K\"ahler forms}

 Recall that a  triple $(\omega, J,
g)$ on a smooth manifold $X$ is called a K\"ahler structure if
$\omega$ is a symplectic form, $J$ is a complex structure, $g$ is a
Riemannian metric such that $g(u,v)=\omega(u, Jv).$ There are
$4-$manifolds admitting complex structures but not admitting
K\"ahler structures, e.g. $S^3\times S^1$. A deep result is that any
complex surface admits a K\"ahler structure if and only if  the
first Betti number is even (see \cite{BPV}). For the analogous
question when  a 4-manifold admitting  a symplectic structure also
admits a K\"ahler structure,
it  has a negative
answer for $\kappa=1,2$: There are many symplectic 4-manifolds with even $b_1$
(or $b_1=0$) admitting no Kaehler structure (see \cite{Go1}).
It does  have a positive answer in the case
$\kappa=-\infty$ (see \cite{Liu1}). We might  be able to answer this
question completely in the case of $\kappa=0$ as well (see
\cite{L1}). This question can also be approached in terms of
Lefschetz fibration structures. The existence of such a structure
essentially characterizes symplectic structures in dimension 4 (see
\cite{Do2}, \cite{Go2}). It is shown in \cite{ST} that all genus 2
Lefschetz fibrations with only irreducible singular fibers
\footnote{the irreducibility condition is necessary by \cite{OS}.}
and transitive monodromy admit K\"ahler structures. A recent result
in \cite{Au} says that, after stably summing sufficiently number of
holomorphic Lefschetz fibrations,
 all  Lefschetz fibrations admit K\"ahler structures.

From now on in this   section,  $X$ is a closed, smooth, oriented
$4-$manifold admitting  a K\"ahler structure. We will compare the
spaces of symplectic structures and K\"ahler structures on $X$.

\subsection{Towards a symplectic Nakai-Moishezon criterion}
Let ${\mathcal C}_{J}$ be the K\"ahler cone of the complex surface $(X,J)$
and $K_J$ be the canonical class.
We will call a complex surface $(X,J)$ a K\"ahler surface if ${\mathcal C}_{J}$ is non-empty.
 Let  $H_J^{1,1}$ denote the real part of the $(1,1)-$subspace of $H^2(X;{\mathbb
C})$ determined by $J$. Then, by the Hodge index theorem, the cone
of classes with positive squares in $H_J^{1,1}$ has two connected
components.  The K\"ahler cone ${\mathcal C}_{J}$ is inside one such
component, which we call the forward cone of $H_J^{1,1}$ and denote
it by  $FH_J^{1,1}$. The K\"ahler Nakai-Moishezon criterion in
\cite{Bu} asserts that, for a complex surface $(X,J)$ \footnote{The
higher dimensional extension is given in \cite{DP}.} with ${\mathcal
C}_J\ne \emptyset$,
$${\mathcal C}_J=\{e\in FH_{J}^{1,1}| \hbox{ $e(C)>0$ for any holomorphic curve
$C$ with $C\cdot C<0$}\}.$$

Thus each holomorphic curve $C$ with $C\cdot C<0$ determines a face of the K\"ahler cone.
 Since such a curve in a minimal
K\"ahler surface would disappear after deforming the complex
structure to a generic nearby almost complex structure,  the
following question was  raised in \cite{Bi1} and \cite{Bi3}.

\begin{question} \label{oneone} For a minimal K\"ahler surface $(X,J)$, is every class of
positive square in $H^{1,1}_{J}$ realized by a symplectic form?
\end{question}

We have recently made some progress towards Question \ref{oneone} in
\cite{LU}. More precisely, the following result is proved using the
normal connected sum construction.

\begin{theorem} \label{negative}
Suppose $C$ is an embedded connected symplectic surface representing
a homology class $E$ with $E\cdot E=-k<0$ and $a=\omega(E)$. Let
$h=k$ if $C$ has positive genus or $C$ is a sphere with $k$ even,
and $h=k+1$ if $C$ is a sphere with $k$ odd. Then there are
symplectic forms $\omega_t$ representing the classes $[\omega]
+tPD(E)$ for any $t\in [0, {2a\over h})$.
\end{theorem}

 To apply
such a construction to Question \ref{oneone} we need to turn
holomorphic curves into embedded symplectic surfaces.  It is shown
in \cite{Mc4} that any irreducible simple pseudo-holomorphic curve
can be perturbed to a pseudo-holomorphic immersion, possibly after a
$C^1$-small change in the almost complex structure. We show how to
further perturb such an immersion to an embedding, which is
$J'$-holomorphic for an almost complex structure arbitrarily
$C^1$-close to $J$. Consequently, we have

\begin{theorem} Let $(X, J)$ be a minimal K\"ahler surface. Then, inside the symplectic cone,
the K\"ahler cone can be enlarged across  any of its open face
determined by an irreducible curve with negative self-intersection.
In fact, if the curve is not a rational curve with odd
self-intersection, then the reflection of the K\"ahler cone along
the corresponding face is in the symplectic cone.
\end{theorem}

In addition, for a minimal surface of general type, the canonical
class $K_J$ is shown to be in the symplectic cone in \cite{Sm2},
\cite{Ca} (It is in the K\"ahler cone ${\mathcal C}_J$ if and only
if there are no $-2$ curves).

\subsection{The case $p_g=0$}
 In this case  we have $b^+=1$ and $H^{1,1}_{J}=H^2(X;{\mathbb R})$.
We have completely determined the symplectic cone in Theorems
\ref{positivecone} and \ref{symplecticcone}. Together with Theorem
\ref{defiso}, we have

\begin{theorem} When $p_g=0$,
Question \ref{oneone} has an affirmative answer. In addition,
 a symplectic form deformation equivalent
to a K\"ahler form in the same cohomology class is itself K\"ahler.
\end{theorem}

\begin{definition} A complex structure $J$ is called symplectic
generic if the K\"ahler cone ${\mathcal C}_{J}$ coincides with the
$K_J-$symplectic cone ${\mathcal C}_{(X,K_J)}$.
\end{definition}

The central question here is

\begin{question} \label{generic}
Does there exist a symplectic generic complex structure?
\end{question}

For manifolds with $\kappa=-\infty$, by Theorem \ref{modulispace}, a
positive answer to this question would imply that every symplectic
form is K\"ahler in this case. It is also worth mentioning that, as
first observed in \cite{McP}, a positive answer for all rational
manifolds would imply the following longstanding Nagata conjecture
(\cite {N}) on the minimal possible degree of a  plane curve with
prescribed singularities.

 \begin{conjecture} Let $p_1,..., p_l\in \mathbb {CP}^2$ be $l\geq
 9$ very general points. Then for every holomorphic curve $C\subset
 \mathbb{CP}^2$ the following inequality holds:
 $$\deg(C)\geq \frac { {\rm mult}_{p_1}C+...+{\rm mult}_{p_k}C}
 {\sqrt k}.$$
 \end{conjecture}

If we let $h$ denote the hyperplane class of $\mathbb{CP}^2$ and
$e_1,..., e_k$ denote the Poincar\'e dual to the the exceptional
divisors. Then, via the K\"ahler Nakai-Moishezon criterion, the
Nagata conjecture is the same as saying that the class
$$\alpha=h-{1\over \sqrt k}(e_1+...+e_k)$$
with $\alpha^2=0$ is in the closure of the K\"ahler cone for a
generic blowup (See \cite{McP}, \cite{Bi6} and \cite{Bi1} for more
extensive discussions). Notice that $\alpha$ is in the closure of
the symplectic cone  by Theorem \ref{packing}.

\begin{prop}
 Symplectic generic complex structures  exist for minimal ruled manifolds and rational manifolds
 with $b^-\leq 9$.
\end{prop}

 For a minimal ruled manifold other than $S^2\times S^2$, any complex
 structure $J$
 arising from a stable rank 2 bundle is symplectic generic (see \cite{Mc5}).
 For rational manifolds with $b^-\leq 8$, any Fano
complex structure is symplectic generic (see \cite{R2}).
 More generally, for  rational manifolds with $b^-\leq 9$,
consider the notion of good generic  surfaces (see e.g. \cite{FM1}).
A good generic surface $(X,J)$ is an algebraic surface such that the
anti-canonical divisor $-K_J$ is effective and smooth, and that any
smooth rational curve has self-intersection no less than $-1$. All
such surfaces are rational surfaces, and for each positive integer
$l$,
 $\mathbb {CP}^2\#l\overline{\mathbb {CP}}^2$  admits
such a structure. For such a surface,  the K\"ahler cone is a nice
subcone of the $K_{J}-$symplectic cone (see \cite{LiL1}),
$${\mathcal C}_J= \{e\in {\mathcal C}_{(X,K_J)}|e\cdot (-K_J)>0\}.$$
When $l\leq 9$, $K_J\cdot K_J\geq 0$, therefore, the condition
$e\cdot (-K_J)>0$ is automatically satisfied for any $e$ in the
$K_J-$symplectic cone.

For related results on K\"ahler $3-$folds with $p_g=0$ see
\cite{R2}, \cite{R3} and \cite{Wil}.


\subsection{The case $p_g>0$}
If $X$ underlies a  minimal K\"ahler surface  with $p_g>0$, then all
complex structures $J$ on $X$ give rise to the same set $\{K_J,
-K_J\}$ by \cite{Wit}. For our purpose  we will denote this set by
$\{K_X, -K_X\}$.

 Let ${\mathcal P}^0={\mathcal
P}$ be the cone of classes of positive squares,  and ${\mathcal
P}^{\alpha}=\{e\in {\mathcal P}|e\cdot \alpha> 0 \}$ for $\alpha\ne
0\in H^2(X;{\mathbb Z})$. The following question, raised in
\cite{LL3}, concerns the (full) symplectic cone of a K\"ahler
surface.

\begin{question} \label{full} If $X$ underlies  a minimal K\"ahler surface
 with $p_g>0$,  ny Theorem \ref{constraint}, the symplectic cone ${\mathcal C}_X$ is
contained in ${\mathcal P}^{K_X}\cup {\mathcal P}^{-K_X}$. Is
${\mathcal C}_X$ equal to ${\mathcal P}^{K_X}\cup {\mathcal
P}^{-K_X}$?
\end{question}

Notice that this question also makes sense in the case $p_g=0$.
Since $K_X^2\geq 0$, ${\mathcal P}^{K_X}\cup {\mathcal P}^{-K_X}$
coincides with ${\mathcal P}$ by the light cone lemma. In
particular, the question has an affirmative answer in this case.
This question is shown to have a positive answer for the product of
the torus and a positive genus Riemann surface in \cite{DL} via the normal connected sum construction.

A related question, raised in \cite{Dr}, is to compare the
cohomology K\"ahler cone with the symplectic cone, where the
cohomology K\"ahler cone is defined to be the union of K\"ahler
cones over all complex structures. There is  a nice answer to both
these questions in the case of K\"ahler surfaces of Kodaira
dimension zero.

 \begin{prop} If $X$ underlies a minimal K\"ahler surface with $K_X=-K_X=0$ and $p_g\geq 1$,
 then  ${\mathcal C}_X={\mathcal P}$. Moreover, the cohomology K\"ahler cone is the same as the symplectic
 cone,
 i.e. every symplectic form is cohomologous to a K\"ahler form.
 \end{prop}

 \begin{proof} It suffices to show that the cohomology K\"ahler cone is
 the positive cone.

By the Kodaira classification of complex surfaces (see e.g.
\cite{BPV}), $X$ is either the 4-torus $T^4$ or the K3 surface, and
each has $b^+=3$. One nice feature for such a manifold is the
existence of a hyperk\"ahler metric $g$. Such a metric induces a
family of complex structures parameterized by the unit sphere of the
imaginary quaternions and the corresponding family of K\"ahler
forms. Denote this sphere by $S^2(g)$, and for each $u\in S^2(g)$,
denote the corresponding K\"ahler form by $\omega_u$. The span $F$
of the $\omega_u$ is a $3-$dimensional positive-definite subspace of
$H^2(X;{\mathbb R})$ (with a basis given by $\{\omega_I, \omega_J,
\omega_K\}$). Let $F^{\perp}$ be the orthogonal complement of $F$,
then $H^2(X;{\mathbb R})=F\oplus F^{\perp}$ as $b^+=3$. The fact we
need now is that $H^{1,1}_{u}={\mathbb R}[\omega_u]\oplus
F^{\perp}$.
 In particular, each class of positive square $e$ is in $H^{1,1}_{u}$ for
some $u\in S^2(g)$. By possibly switching $u$ to $-u$, we can assume
that $e$ is in  $FH^{1,1}_{u}$. In the case of $T^4$, since there
are no holomorphic curves of negative self-intersections for any
complex structure, we conclude from the K\"ahler Nakai-Moishezon
criterion that $e$ is in the K\"ahler cone for $u$.

The same argument almost works for the K3 surface, except that there
are rational curves with self-intersection $-2$ for some $u$. So one
$S^2-$family is not enough, and we will need to use all such
$S^2-$families. According to the surjectivity of the refined period
map (see e.g. Theorem 14.1 in \cite{BPV}), we have

\begin{lemma} \label{period} $e\in{\mathcal P}$ is a K\"ahler class if and only if
there is a positive-definite $2-$plane $U$ in $H^2(X;{\mathbb R})$
such that $e\perp U$, and $e\cdot d\ne 0$ for any integral class $d$
in $H^2(X:{\mathbb Z})$ with $d^2=-2$ and $d\perp U$.
\end{lemma}

Let Gr$^+_e$ be the Grassmannian of positive-definite $2-$planes
which are orthogonal to $e$. Let $\Delta_e$ be the set of $d$ with
$e\cdot d=0$. For  $d\in \Delta_e$, the $2-$planes in Gr$^+_e$ which
are orthogonal to $d$ is a sub-Grassmannian.  Since the complement
of the countable union of these sub-Grassmannians over $\Delta_e$ is
non-empty,  $e$ is K\"ahler by Lemma \ref{period}.
\end{proof}


There has long been a speculation that any symplectic form is still
K\"ahler in this case. Moreover, since the moduli space of complex
structures is connected, the space of K\"ahler forms is also
connected.

 If every class of positive square of $X$
 is in a $H_{J}^{1,1}$ subspace for a complex structure $J$, then the cohomology K\"ahler cone
 of $X$
 contains an open subset of the symplectic cone ${\mathcal C}_X$. The minimal Elliptic surfaces $E(n)$
are likely to have this property.
 For  a surface with such a property, if Question \ref{oneone} can be answered positively,
then its cohomology K\"ahler cone agrees with  ${\mathcal
P}^{K_X}\cup {\mathcal P}^{-K_X}$. Consequently, for such a surface, Question \ref{full} has a positive answer
and every symplectic form is cohomologous to a K\"ahler form.

For  minimal K\"ahler surfaces of general type with $p_g\geq 1$, it
is observed in \cite{Dr} that  the cohomology K\"ahler cone is
strictly smaller than the symplectic cone. On the one hand, by the
Hodge index theorem,
 the cohomology K\"ahler cone of $X$ is contained in the
cone
$$
\{e\in {\mathcal P}|(e^2)(K_X^2)\leq (e\cdot K_X)^2\} \cup {\mathbb R}K_X.
$$
On the other hand, for  a given complex structure $J$, if $\omega$
is a K\"ahler form and $\phi$ is a holomorphic $2-$form, then
$\tau_s=\omega+s{\rm Re}\phi$ is a symplectic form for any $s>0$ due
to the pointwise orthogonality of $(1,1)-$forms and $(2,0)-$forms.
But if $s$ is large enough, then
$$([\tau_s]\cdot K_X)^2=([\omega]\cdot K_X)^2<([\tau_s]^2)(
K_X)^2=([\omega]^2+s^2[{\rm Re}(\phi)]^2)(K_X^2),$$ because
$K_X^2>0$. Therefore, for $s$ large,  the class $[\tau_s]$ is in the
symplectic cone, but not in the cohomology K\"ahler cone.

Nevertheless this `negative' observation can be used to construct
examples of
 many non-holomorphic Lefschetz pencils
on K\"ahler surfaces as follows. Consider a ball quotient $X$ with
$p_g>0$.
 It has a unique complex structure (hence the cohomology K\"ahler cone is of lower dimension
 and does not even contain
 an open subset of the symplectic cone).
   Pick an
integral symplectic structure whose class is outside the $H_J^{1,1}$
subspace for the unique complex structure $J$. According to
\cite{Do2}, there is a Lefschetz pencil whose members are Poincar\'e
dual to large integral multiples of $[\omega]$. Such a Lefschetz
pencil is not holomorphic since  the members of a holomorphic
pencil are holomorphic curves.

\subsection{Donaldson's program}
Recently Donaldson asks the following question in \cite{Do3},

\begin{question}
Suppose $X$ is a compact K\"ahler surface with K\"ahler form $\omega_0$. If $\omega$ is any other symplectic form on $X$, with the same Chern class
and with $[\omega]=[\omega_0]$, is there a sdiffeomorphism $f$ of $X$ with
$f^*(\omega)=\omega_0$?
\end{question}

Donaldson outlines a line of attack involving an almost K\"ahler analogue
of the famous complex  Monge-Amp\`ere equation solved by Yau.
Partial existence result is obtained in \cite{We}. This program seems to be very exciting and ties
up with many aspects discussed in this survey.


\end{document}